\documentclass{article}
\usepackage{xcolor}
\usepackage{transparent}
\usepackage{amsmath,enumitem}
\usepackage{amsthm}
\usepackage[pdftex]{graphicx}
\usepackage{caption}
\usepackage{subcaption}
\usepackage{amssymb}
\usepackage{rotating}
\usepackage{color,soul}
\usepackage{hyperref}
\usepackage{amsmath}
\usepackage{longtable}
\usepackage{tabularx,ragged2e,booktabs}
\usepackage{authblk}
\usepackage[a4paper, margin=3.5cm]{geometry}
\usepackage{floatrow}

\restylefloat{figure}
\restylefloat{table}

\newcommand{\head}[1]{\textnormal{\textbf{#1}}}
\newcolumntype{s}{>{\hsize=.5\hsize}X}

\title{Best Practices for Comparing Optimization Algorithms}
\author[1]{Vahid Beiranvand}
\author[2]{Warren Hare}
\author[1]{Yves Lucet \footnote{corresponding author: yves.lucet@ubc.ca}}
\affil[1]{\small Computer Science, University of British Columbia, Kelowna, BC, Canada}
\affil[2]{Mathematics, University of British Columbia, Kelowna, BC, Canada}
\date{\today}

\begin {document}
\maketitle
\begin{abstract}
The final publication is available at Springer via http://dx.doi.org/10.1007/s11081-017-9366-1

\noindent Comparing, or benchmarking, of optimization algorithms is a complicated task that involves many subtle considerations to yield a fair and unbiased evaluation. In this paper, we systematically review the benchmarking process of optimization algorithms, and discuss the challenges of fair comparison. We provide suggestions for each step of the comparison process and highlight the pitfalls to avoid when evaluating the performance of optimization algorithms. We also discuss various methods of reporting the benchmarking results. Finally, some suggestions for future research are presented to improve the current benchmarking process.

\smallskip
\noindent \textbf{Keywords.} Benchmarking; Algorithm comparison; Guidelines; Performance; Software; Testing; Metric; Timing; Optimization algorithms.
\end{abstract}

\section{Introduction}

As the number of optimization methods, and implementations of those methods, has increased, researchers have pursued comparative studies to evaluate their performance.  When done well, such studies can be of great value in helping end-users choose the most suitable optimization method for their problems.  Such studies are generally referred to as optimization benchmarking.

In the most general sense, benchmarking is the comparison of one or more products to an industrial standard product over a series of performance metrics.  In the case of benchmarking optimization algorithms, the products are the specific implementations of given algorithms, and the performance metrics are generated by running the implementations on a series of test problems.  This framework presents a certain clarity in benchmarking optimization algorithms, as there is at least some agreement on what constitutes ``better''.  If one algorithm runs faster, uses less memory, and returns a better final function value, on all possible problems, then it can be considered better than the alternative.  Of course, in practice such a clear conclusion seldom arises. Thus, interpreting the conclusions of algorithmic comparisons can be tricky.

Nonetheless, when done well, benchmarking optimization algorithms can have great practical value.  It can reveal both strengths and weaknesses of an algorithm, which allows for better research focus.  It can aid in determining if a new version of optimization software is performing up to expectations.  And, it can help guide end-users in selecting a good choice of algorithm for a particular real-world problem.

However, when done poorly, benchmarking optimization algorithms can also be misleading.  It can hide algorithm's weaknesses (or strengths), report improvements that do not exist, or suggest the incorrect algorithmic choice for a given situation.

In optimization benchmarking many subjective choices are made, such as the test set to solve, the computing environment to use, the performance criteria to measure, etc. Our primary objective in this paper is to help researchers to design a proper benchmarking approach that is more comprehensive, less biased, and less subject to variations within a particular software or hardware environment. Our secondary objective is to provide a comprehensive review of the benchmarking literature for optimization algorithms.

In pursuing these objectives, we focus on single-objective optimization algorithms that run in serial (i.e., that do not use parallel processing).  Comparing algorithms for multi-objective optimization, or optimization algorithms that use parallel processing, introduce new levels of complexity to the benchmarking process.  While we provide some comments on the challenges for benchmarking some algorithms in the conclusions, we consider these issues outside of the scope of this paper.

We also note that much of the presentation within this paper discusses algorithms as if the underlying optimization problem is a continuous unconstrained problem.  This is for ease of presentation, and in most cases translating the ideas to other styles of optimization problems is clear.  As such, we limit ourselves to discussing specific styles of optimization problems only when the translation is not straight-forward.

\subsection{Historical overview of benchmarking in optimization}

We begin with a brief historical overview of optimization benchmarking.

One of the very first studies in benchmarking of algorithms was given by Hoffman et al., in 1953 \cite{117}, in which three different methods for linear programming were compared. Although this computational experiment was performed early in the development of computers, when there existed almost no compiler and programming environment, the reported techniques have been used for a long time and can be considered as the foundation of the current comparison techniques. They include such ideas as, using test sets to compare algorithms, employing performance measures (accuracy, CPU time, number of iterations, and convergence rate), and paying attention to the impact of coding on the performance of algorithms.

Another early contribution to the field of benchmarking is Box's work from 1966 \cite{26}.  In this work, Box evaluates the performances of eight algorithms for unconstrained optimization using a collection of 5 test problems with up to 20 variables. He considered the number of function evaluations, the importance of model size and the generality of the optimization algorithms.

In the late 1960's, optimization benchmarking research began to expand rapidly. Comparative studies have been performed throughout the optimization literature, for example in unconstrained optimization \cite{59, 53, 23}, constrained optimization \cite{60, 56, FamularoPuglieseSergeyev2002}, nonlinear least squares \cite{54, 55, 19, 109}, linear programming \cite{104, 145, 47}, nonlinear programming \cite{103, 54, 62, 51, 106, 107, 108, 105, 18, 19, 17, 21, 75, 102, 70}, geometric programming \cite{64, 58}, global optimization \cite{88, 78, 123, 80, 73, 22, SS00, ZZ08, KM16}, derivative-free optimization \cite{82, 95, 99, SK06, PSKZ14}, and other areas of optimization \cite{59, 53, 61, 120,  11,148, 68} -- amongst many more.

In addition, a few researchers have focused on improving the (optimization) benchmarking process.  In 1979, Crowder et al.~\cite{38}, presented the first study that attempted to provide standards and guidelines on how to benchmark mathematical algorithms.  It includes a detailed discussion of experimental design and noted the necessity of {\em a priori} experimental design. The authors paid attention to reproducibility of the results and provided a method for reporting the results. In 1990, similar research conducted by Jackson et al.\ \cite{9} delivered an updated set of guidelines.  In 2002, Dolan and Mor\'e introduced performance profiles \cite{2}, which have rapidly become a gold standard in benchmarking of optimization algorithms with more recent work pointing out its limitations~\cite{GOULD-16}.  In this paper, we attempt to provide a modern picture of best-practice in the optimization benchmarking process.

\subsection{Paper framework}

We now provide a general framework for benchmarking optimization algorithms, which we also use to structure discussions in the paper.
\begin{enumerate}
\item \textbf{Clarify the reason for benchmarking.} In Section \ref{sec:reason}, we discuss some of the common reasons to compare optimization algorithms, and some of the pitfalls that arise when the purpose of benchmarking is unclear.
\item \textbf{Select the test set.} In Section \ref{sec2}, a review of test sets for various problem categories is presented, the challenges related to test sets are discussed, and some guidelines are provided for assembling an appropriate test set.
\item \textbf{Perform the experiments}.  In Section \ref{sec3}, we review and discuss various considerations related to the critical task of designing experiments, including performance measures, tuning parameters, repeatability of the experiments, and ensuring comparable computational environments.
\item \textbf{Analyze and report the results}. Section \ref{sec4} contains a review of different reporting methods for optimization algorithms, including tabular methods, trajectory plots, and ratio-based plots (such as performance and data profiles).
\end{enumerate}

In addition to the aforementioned sections, Section \ref{sec5} contains a review of recent advances in the field of automated benchmarking and Section \ref{sec6} presents some concluding thoughts.

\section{Reason for benchmarking}\label{sec:reason}

Having a clear understanding of the purpose of a numerical comparison is a crucial step that guides the rest of the benchmarking process.  While seemingly self-evident, it is surprisingly easy to neglect this step.  Optimization benchmarking has been motivated by a variety of objectives.  For example:
    \begin{enumerate}
        \item To help select the best algorithm for working with a real-world problem.
        \item To show the value of a novel algorithm, when compared to a more classical method.
        \item To compare a new version of optimization software with earlier releases.
        \item To evaluate the performance of an optimization algorithm when different option settings are used.
    \end{enumerate}
In a practical sense, all of these work towards gathering information in order to rank optimization algorithms within a certain context.  However, the context can, and should, play a major role in guiding the rest of the benchmarking process.  

For example, if the goal is to select the best algorithm for a particular real-world application, then the test problems (Section \ref{sec2}) should come from examples of that application.  

Alternately, if the goal is to show the value of a new optimization algorithm, then it is valuable to think about exactly where the algorithm differs from previous methods.  Many new algorithms are actually  improvements on how a classical method deals with some aspect of an optimization problem.  For example, in \cite{c14} the authors develop a new method to deal with nonconvexity when applying a {\em proximal-bundle method} to a {\em nonsmooth optimization problem}.  As such, to see the value of the method, the authors compare it against other proximal-bundle methods on a collection of nonconvex nonsmooth optimization problem.  If they had compared their method against a quasi-Newton method on smooth convex optimization problems, then very little insight would have been gained.

Regardless of the reason, another question researchers must consider is what aspect of the algorithm is most important.  Is a fast algorithm that returns infeasible solutions acceptable?  Is it more important that an algorithms solves every problem, or that its average performance is very good?  Is the goal to find a global minimizer, or a highly accurate local minimizer?  The answers to these questions should guide the choice of performance metrics that need to be collected (Section \ref{sec3}) and how they should be analyzed (Section \ref{sec4}).  Answering these types of questions before running the experiments is time well spent.

\section{Test sets} \label{sec2}
%suggestions on how to build a "good" test set and links to available sets
A test problem contains a test function together with some further criteria such as the  constraint set, feasible domain, starting points, etc. A test set is a collection of test problems. Obviously, benchmarking only yields meaningful results when competing algorithms are evaluated on the same test set with the same performance measures.

The selection of the appropriate test sets to benchmark the performance of optimization algorithms is a widely noticed issue among researchers \cite{72, 8, 9, SSL13, ZZ08}. Generally, there are three sources for collecting test problems: real-world problems, pre-generated problems, and randomly-generated problems. Real-world problems can be found through instances of specific applications, and pre-generated problems exist in common test set libraries, see Table \ref{Table2}.  Conversely, randomly-generated test problems are often more {\em ad hoc} in nature, with researchers creating methods to randomly generate test problems that are only used in a single paper \cite{13, c14, c19} (among many others).  However, some researchers have gone to the effort to study methods to randomly-generate test problems for a given area; some examples appear in Table \ref{Table2b}.

While the real-world test sets provide specialized information about the performance of the optimization algorithms within a specific application, the results may be difficult to generalize.  The difficulties lie in the facts that real-world test sets are often small and the problems are often application-specific.   Nonetheless, if the goal is to determine the best algorithm to use for a particular real-world application, then a real-world test set focused on that application is usually the best option.

On the other hand, the artificial and randomly-generated test sets can provide useful information about the algorithmic characteristics of optimization algorithms.  Artificial and randomly-generated test sets can be extremely large in size, thereby providing an enormous amount of comparative data.  However, it can be difficult to rationalize their connection to the real-world performance of optimization algorithms.  If the goal is to compare a collection of algorithms across a very wide spectrum, then artificial and randomly-generated test sets are usually the better option.

When selecting a test set, it is always important to keep the particular goal of the comparison in mind.  Regardless of the goal, an appropriate test set should generally seek to avoid the following deficiencies.

\begin{enumerate}
\item[i.] \textbf{Too few problems.} Having more problems in the test set makes the experiment more reliable and helps the results reveal more information about the strengths or weaknesses of the evaluated algorithms.

\item[ii.] \textbf{Too little variety in problem difficulty.} A test set containing only simple problems is not enough to identify the strengths and weaknesses of algorithms. In contrast, a test set which only has problems that are so difficult that no algorithm can solve them, clearly, does not provide useful information on the relative performance of algorithms.

\item[iii.] \textbf{Problems with unknown solutions.} When possible, it is better to use test problems where the solution is known.  Depending on the analysis performed (see section \ref{sec4}), the ``solution'' could be interpreted as the minimum function value, or the set of global minimizers.   Having access to the solution greatly improves the ability to evaluate the quality of the algorithmic output.  However, when the test set is comprised of real-world test problems, then a lack of known solutions may need to be accepted as inevitable.

\item[iv.] \textbf{Biased starting points.} Allowing different algorithms to use different starting-points will obviously bias the result.  However, more subtle problems can also exist.  For example, if a starting point lies on the boundary of a constraint set, then an interior point method will be severely disadvantaged.  Another example comes from considering the Beale test function, which has a global minimizer at $(3,0.5)$ \cite{23}.  If a {\em compass search} (see, e.g.,  \cite{KoldaLewisTorczon2003}) with an initial step length of $1$ is started at $(0.5, 0.5)$, then it will converges to the exact minimizer in just 4 iterations.  However, if a starting point of $(0.51, 0.51)$ is used, then the exact same algorithm will use $63$ iterations to reach a point within $10^{-2}$ of the global minimizer.\footnote{Note that this example is artificially constructed to emphasize the results; the recommended starting point for the Beale test problem is $(1,1)$.}

\item[v.] \textbf{Hidden structures.}  Many test sets have some structure that is not realistic in real-world problems.  For example, about 50\% of the problems in the test set \cite{23} have solution points that occur at integer valued coordinates.  An algorithm that employs some form of rounding may perform better than usual on these problems.
\end{enumerate}

Thus, when choosing test sets for the benchmarking task the following considerations should be taken into account as much as possible.
\begin{enumerate}
\item[i.] If the test set contains only few problems, then the experiment should be referred to as a \textit{case study} or a \textit{proof of concept}, but not benchmarking.  While there is no fixed number that determines what is ``enough problems to be considered benchmarking'', we recommend that in order to achieve a reliable conclusion about the performance, an experiment should contain at least 20 test problems (preferably more). In the specific case of comparing a new version of an optimization algorithm with a previous version, the number of test problems should be significantly greater -- in the order of 100 or more.  In all cases, the more problems tested, the more reliable the conclusions.

\item[ii.] When possible, a test set should include at least two groups of problems: an \textit{easy group} which consists of the problems that are easy to solve within a reasonable time on a regular contemporary computer using all the optimization algorithms tested, and a \textit{hard group} that contains the problems which are solvable but computationally expensive and may require a specific optimization algorithm.

\item[iii.]  Whenever possible, ensure at least a portion of the test set includes problems with known solutions.

\item[iv.]  For test sets that include starting points, new starting points can be generated by introducing a small (random) perturbations to the given starting points.  For other test sets, randomly-generated starting points can be created from scratch.  In either case, all starting points should be created for each problem, and then every algorithm should be provided the same starting point for testing.  This approach can be further used to increase result reliability, by repeating tests on the same function with a variety of starting points.

\item[v.]  Examine the test set with a critical eye and try to determine any hidden structure.  Some structures can be removed through methods similar to the random perturbation of starting points in (iv).  One quick test is to set an algorithm to minimize $f(x)$ starting at $x^0$ and then set the algorithm to minimize $\hat{f}(x) = f(x - p)$ starting from $\hat{x}^0 = x^0 - p$ (where $p$ is any random point).  Constraints can then be shifted in a similar manner, effectively shifting the entire problem horizontally by the vector $p$.  While it relocates the origin, and moves any constraints away from special integer values, it has no theoretical effect on the geometry of the problem.  As such, the results of both tests should be extremely similar (theoretically they should be identical, but numerical errors may cause some deviation).  If the results of both tests differ, then perhaps some hidden structure is being exploited by the algorithm, or perhaps some hidden constraints are causing issues.  Regardless of the reason, the researcher should recognize the issue and consider a wider test set.
\end{enumerate}

Using suitable standard test sets is usually a good option when benchmarking optimization algorithms.  In particular, it is usually easier to compare results across research groups when standard tests are employed, although even within a specific research field there is generally no consensus on the appropriate test set to draw specific conclusions. Many interesting and diverse test sets have been reported in the literature, see Tables \ref{Table2} and \ref{Table2b}.

\begin {table}[H]
\captionof{table}{Some test sets reported in the literature.\label{Table2}}
\begin{center}
\begin{tabularx} {\textwidth}{ss}
  \toprule[1.5pt]
	\head{Collection type} & \head{Resources} \\
  \toprule[1.5pt]
  Unconstrained optimization problems & \cite{4, 23, 63, 118} \\   \midrule
  Global optimization & GAMS \cite{72}, COCONUT \cite{66}, and other collections 
  \cite{98, 90, 89, 65, 94, 96, FamularoPuglieseSergeyev2002} \\   \midrule
  Linear programming & \cite{c20} \\ \midrule
  Local optimization &\cite{90, 98} \\ \midrule
  Nonlinear optimization problems & CUTEr \cite{85,86}, CUTEst \cite{111}, COPS \cite{6,8}, and  collections \cite{52, 83, 84, 57, 142} \\ \midrule
  Mixed integer linear programming & MIPLIB \cite{PEKO2002, 49} \\
  \bottomrule[1.5pt]
\end{tabularx}
\end{center}
\end{table}

Producing random test sets using test problem generators has its own drawbacks. Are the generated problems representative or difficult? Is there any hidden structure in the problems? Some papers that use random test problem generators are listed in Table~\ref{Table2b}.

\begin {table}[H]
\captionof{table}{Some test problem generators reported in the literature.\label{Table2b}}
\begin{center}
\begin{tabularx} {\textwidth}{ss}
  \toprule[1.5pt]
  	\head{Test problem generators} &  \\  	
  \toprule[1.5pt]
    Network programming & \cite{157}\\\midrule
    Nonlinear optimization problems & \cite{51,91}  \\ \midrule
    Combinatorial problems & \cite{131, 157} \\
    \midrule
    Quadratic Programming & \cite{156} \\
    \midrule
    Global Optimization & \cite{89, 93, 94, NL14} \\
  \bottomrule[1.5pt]
\end{tabularx}
\end{center}
\end{table}

\begin{figure}[ht]
\centering
  \def\svgwidth{\columnwidth}
	\footnotesize %\small
  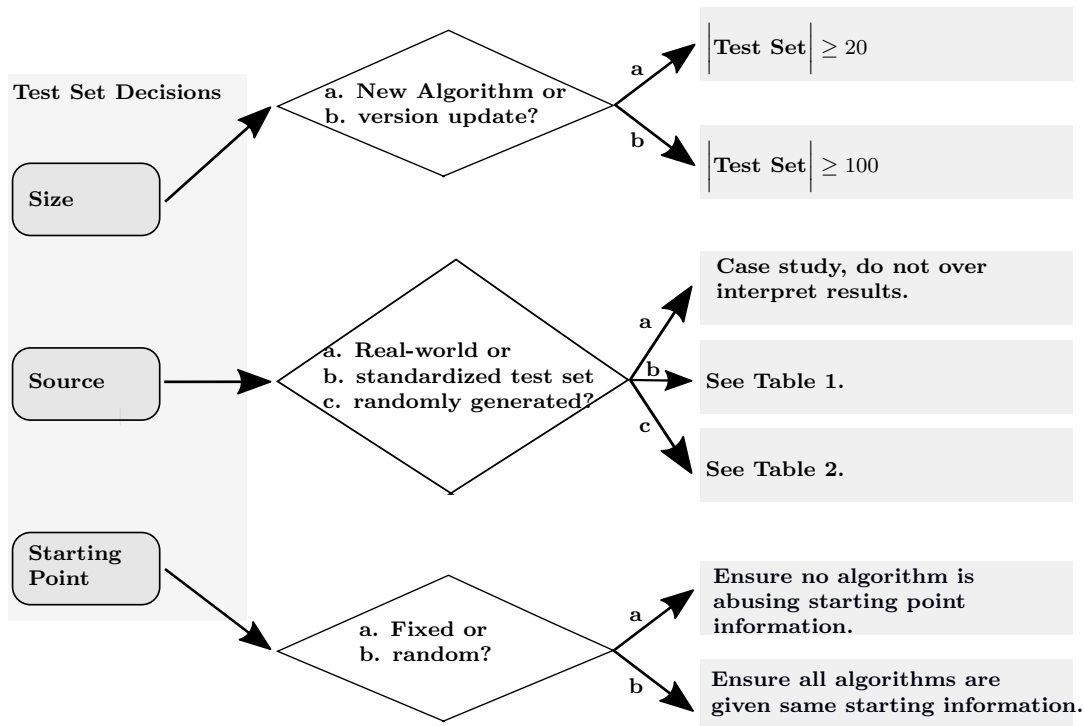
  \caption{Test set decision tree.}
  \label{Fig06TestSet}
\end{figure}

Figure \ref{Fig06TestSet} shows a decision tree that summarizes the fundamental decisions required for assembling an appropriate test set for benchmarking of optimization algorithms.

\section{Performing the experiments}\label{sec3}

The performance of algorithms is influenced by two general types of factors: \textit{environmental factors} and \textit{algorithmic factors}.

\textit{Environmental factors} refer to factors that are out of the algorithm scope and usually beyond the control of the researcher. A common example is the computer environment used to test the algorithms, which includes processor speed, operating system, computer memory, etc. Environmental factors may also include the programmer's skill and the programming language/compiler used. This is particularly evident when multiple pieces of software by a variety of programmers are being compared.  In essence, if the benchmarking process is repeated by another researcher elsewhere, then environmental factors are unlikely to remain constant, and so the benchmarking results are expected to change.

\textit{Algorithmic factors} are related to the algorithm itself. These are factors that are considered global across a variety of computing platforms. If the software is programmed by the researcher, then it is assumed these factors are independent from the implementation aspects of the algorithm.

\emph{Optimization benchmarking} seeks to measure the algorithmic factors, and proceeds under the key \emph{assumption} that, while environmental factors are expected to change the results, the algorithmic factors are sufficiently strong that the general ranking of algorithms should remain constant under the specific ranges of parameters under consideration.

To compare algorithms, it is necessary to collect data that measures the overall performance of each algorithm.  This is done by running each algorithm on the test set (discussed in Section \ref{sec2}), and collecting data on the results.  The data collection and the selection of performance measures is based on the research questions motivating the experimental study.  In general, performance measures fall into 3 categories: efficiency, reliability, and quality of algorithmic output. We discuss these performance categories in Subsections \ref{ss:efficiency}, \ref{ss:reliability}, and \ref{ss:qualityofsolution}.  Table \ref{Table3} provides a classification of the comparative measures for optimization algorithms based partly on the guidelines provided by Hoffman and Jackson \cite{124}.

\begin {table}[H]
\caption{Comparative measures.}
\centering
\begin{tabularx} {\textwidth}{lX}
  \toprule[1.5pt]
	\head{Performance Category} & \head{Example Criteria} \\
  \toprule[1.5pt]
        & 1. Number of fundamental evaluations \\
  Efficiency & 2. Running time \\
        & 3. Memory usage \\
        \toprule[1.5pt]
          & 1. Success rate \\
  Reliability & 2. Number of constraint violations \\
        & 3. Percentage of global solutions found \\
        \toprule[1.5pt]
        & 1. Fixed-cost solution result \\
  Quality of Solution & 2. Fixed-target solve time \\
        & 3. Computational accuracy \\
        \toprule[1.5pt]
\end{tabularx}
\label{Table3}
\end{table}

In order to allow for maximal data analysis (and thereby the best understanding of overall performance), it is recommended to collect at least some data from every performance category.

\subsection{Efficiency}\label{ss:efficiency}

The efficiency of an optimization algorithm refers to the computational effort required to obtain a solution. In mathematical programming, there are two primary measures of efficiency, the number of fundamental evaluations and the running time.  A third, less common, measure is memory usage.

\textbf{Number of fundamental evaluations}: The term fundamental evaluation is used to refer to any subroutine that is called by the algorithm in order to gain fundamental information about the optimization problem.  The most obvious example is an objective function evaluation but the evaluation may involve complex simulation algorithms.  Other fundamental evaluations could include gradient evaluations, Hessian evaluations, or constraint function evaluations.  The number of fundamental evaluations can be used as a standard unit of time, and is often assumed to be platform independent. In many situations, the number of fundamental evaluations is a particularly important measure, as for real-world problems these evaluations often dominate over the internal workings of the algorithm 
\cite{27, 38, c3, 7, 25, 13, 36,53, 61, 62, 66, 73, 105, 129, E85, PSKZ14, KS15}.
%\cite{27, 38, 7, 25, 13, 36, 66, 73, 129, PSKZ14, KS15}. 
Note however that this measure is unreasonable when fundamental evaluations do not dominate the internal workings of the algorithm \cite{4}.

\textbf{Running Time}: Running time, as a measure for optimization benchmarking, is usually measured by either CPU time or wall clock time.\footnote{\textit{Wall clock time} refers to the amount of time the human tester has to wait to get an answer from the computer.  Conversely, \textit{CPU time} is the amount of time the CPU spends on the algorithm, at the exclusion of operating system tasks and other processes.}  Wall clock time contains CPU time, and has been argued to be more useful in real-world setting \cite{109}. However, Wall clock time is not reproducible or verifiable since it is tied to a specific hardware platform and software combination. CPU time is considerably more stable, as it is independent of background operations of the computer.  Moreover, CPU time is more-or-less consistent for the same version of an operating systems running on the same computer architectures. %Thus, the majority of optimization benchmarking is based on CPU time, with the premise that this time provides a reasonable starting point to predict algorithmic performance in terms of wall clock time \cite{12, 13, 57, 58, 62, 75, 109, 105, 106}.

It should be noted that, due to the wide variety and complexity of modern computing architectures,  the number of situations in which time is dominated by memory access costs is increasing, hence the precision of CPU timers has been reduced. To improve the precision of CPU timers, tools such as cache and memory access tracers can help obtaining a more accurate CPU time performance. For a more detailed discussion on these techniques we refer to \cite{c1, 112}.

Another issue with CPU time is the increasing prevalence of multi-core machines. Thorough reporting would require indicating the number of cores, the CPU time for each core, but also how efficiently the different levels of memory were used and cache hits/misses. Since such measurements are not straightforward to obtain for multi-core machines, the wall-clock time along with the hardware specifications are usually reported. (Unless the new algorithm contribution focuses specifically on optimizing computation for a multi-core architecture, in which case more precise measures are warranted.) Eventually, the onus is on the researchers to explain how simplified measurements support the conclusions drawn; this is especially true for multi-core machines. 

Regardless of whether wall clock time or CPU time is used, in order to maximize the reliability of the data, it is important to ensure that any background operations of the computer are kept to a minimum.  Furthermore, any manuscript regarding the benchmarking should clearly state which form of running time was collected.

\textbf{Other measures:} In addition to the categorization presented above, in some specific cases, there is another issue that influences the choice of an appropriate measure for running time: \textit{the type of algorithm}. For example, to evaluate the running time for branch-and-bound based algorithms, the number of branch-and-bound nodes is a common criterion, while for simplex and interior point based algorithms, the number of iterations is often used. Therefore, when deciding on the choice of a suitable efficiency measure, the type of algorithms to be evaluated should also be taken into account.

\subsection{Reliability}\label{ss:reliability}

The reliability and robustness of an optimization algorithm is defined as the ability of the algorithm to ``perform well'' over a wide range of optimization problems \cite{23, 81}. The most common performance measure to evaluate the reliability is success rate \cite{c4, 58, 105, SS00}.  Success rate is gauged by counting the number of test problems that are successfully solved within a pre-selected tolerance.  This can be done using objective function value, or distance of the solution point from a minimizer.  In convex optimization these two approaches are largely, but not perfectly, interchangeable.  However, if the objective function has multiple local minimizers, then the researcher must decide whether good local solutions are acceptable outcomes, or if the algorithm must converge to a global minimizer \cite{51, 55}.  In addition to the success rate, the average objective function values and the average constraint violation values have also been reported to measure reliability \cite{51}.

When studying reliability, the researcher should consider whether the algorithms are deterministic, or non-deterministic, and repeat tests multiple times if the algorithm is non-deterministic. Reliability can be based on a fixed starting point (if one is given with the test set), but it is often better to use multiple starting points.

%\textbf{Deterministic versus Non-deterministic:} An algorithm is referred to as {\em deterministic} if, whenever the algorithm is given the same initial inputs (objective function, starting point, any algorithmic parameters), it will produce the same result.  Any algorithm that does not satisfy this is called {\em non-deterministic}.
%Typical examples of non-deterministic methods are {\em genetic algorithms} \cite{c23}, where each iteration involves random {\em mutations} of current candidate solutions.  As the mutations are random, rerunning the algorithm with the same input may result in vastly different outputs. A good rule-of-thumb is, if the algorithm includes a random number at any point, then it is non-deterministic.  Otherwise, it is probably deterministic.

In deterministic optimization algorithms, reliability can be interpreted as the number of problems in the given test set that are solved by the optimization algorithm. When dealing with non-deterministic algorithms, it is important to repeat each test multiple times, to make sure that reliability is measured in aggregate, and not skewed by a single lucky (or unlucky) algorithmic run.  

Using multiple repeats of each test raises the issue of how to aggregate the results.  One option is to consider each algorithmic run as a separate test problem and then compare solvers across this expanded test set.  This allows comparisons based on worst-case or best-case scenarios.  Another option is to use averaged data, for example, average runtime, average solution accuracy, average reliability, etc.  If averaging is used, then it is important to also include standard deviations of the data.  In either case, data collection is best performed by considering each algorithmic run as a separate test problem, as average values can easily be extracted from this data, while reconstructing the full test data from averaged values is not possible.

It should be noted that, in some cases multiple repeats of a non-deterministic method is impractical due to the time it takes to solve a single problem.

\textbf{Multiple starting points:}  As mentioned in Section \ref{sec2}, many academic test problems come with suggested starting points.  While algorithms should always be tested using these starting points, it is often enlightening to test the algorithm using other starting points as well.  Most deterministic algorithms should show little change in performance if a starting point is perturbed by a small random vector -- provided the new starting point retains whatever feasibility properties the algorithm requires in a starting point.

Hillstrom \cite{15} is one of the first to consider testing optimization algorithms at nonstandard starting points.  He proposed using random starting points chosen from a box surrounding the standard starting point.  In another approach to this problem, in \cite{23} the authors present a large collection of test functions along with some procedures and starting points to assess the reliability and robustness of unconstrained optimization algorithms.  In some cases, prior knowledge is available about the solution of a test problem. Some methods use such information to construct a starting point close to the optimal solution \cite{150}.

Regardless of how starting points are selected, fair benchmarking requires all the algorithms to use the same starting point for each test. Therefore, starting points should be generated and stored outside of the testing process.

\subsection{Quality of algorithmic output}\label{ss:qualityofsolution}

The quality of the algorithmic output is obviously important when comparing optimization algorithms.  Measuring quality falls into two easily separated categories: a known solution is available, and no known solutions are available.

\textbf{Known solution available:} When the expected solution for a problem is available, two methods can be employed to measure the quality of an algorithmic output: fixed-target and fixed-cost \cite{1, 41, 95}.

In the \textit{fixed-target} method, the required time (function calls, iterations, etc) to find a solution at an accuracy target $\varepsilon_{target}$ is evaluated. The main problem with fixed-target methods is that some algorithms may fail to solve a test problem. Therefore, the termination criterion cannot rely only on accuracy, but should also include some safety breaks such as the maximum computational budget. If the algorithm successfully reaches the desired accuracy, then the time to achieve the accuracy can be used to measure the quality of the algorithm on that test problem.  If the algorithm terminates before reaching the desired accuracy, then it should be considered unsuccessful on that test problem.

Let $x^0$ be the initial point from a test run, $\bar{x} \in \mathbb{R}^n$ be the termination point obtained from the test run, and $x^* \in \mathbb{R}^n$ be the known solution for the problem.  In the \textit{fixed-cost} approach, the final optimization error $f(\bar{x})-f(x^*)$ is checked after running the algorithm for a certain period of time, number of function calls, number of iterations, or some other fixed measurement of cost.  Then, the smaller the final optimization error is, the better the quality of the algorithmic output.

The fixed-target versus fixed-cost decision can be seen as a multiobjective problem. It is analogous in engineering to minimizing cost, subject to constraints on performance versus maximizing performance, subject to a constraint on cost.

If a fixed-cost approach is used, then there are many options on how to quantify the accuracy of the algorithmic output.  We need to determine whether $\bar{x}$ approximates $x^*$ or not. For example, this can be done using the function value or the distance from the solution:
    $$f_{\tt acc} = f(\bar{x}) - f(x^*), \quad \mbox{and} \quad x_{\tt acc} = \|\bar{x} - x^*\| ~~\mbox{respectively}.$$
It is often valuable to ``normalize'' these quantities by dividing by the starting accuracy:
    $$f_{\tt acc}^n = \frac{ f(\bar{x}) - f(x^*)}{f(x^0) - f(x^*)}, \quad \mbox{and} \quad x_{\tt acc}^n = \frac{\|\bar{x} - x^*\|}{\|x^0 - x^*\|}.$$
Finally, to improve readability, and reduce floating point errors, many researchers take a base-10 logarithm:
    $$\begin{array}{rcl}
    f_{\tt acc}^l &=& \log_{10}(f(\bar{x}) - f(x^*)) - \log_{10}(f(x^0) - f(x^*)), \quad \mbox{and} \\
    x_{\tt acc}^l &=& \log_{10}(\|\bar{x} - x^*\|) - \log_{10}(\|x^0 - x^*\|).\end{array}$$
The values $f_{\tt acc}^l$ and $x_{\tt acc}^l$ can be loosely interpreted as the negative of the number of new digits of accuracy obtained (measured on a continuous scale), thus making these values very useful for discussion. Finally, to avoid exact solutions making an algorithm look better than it is, one can select a ``maximal improvement value'' $M$ (typically about 16) and set
    \begin{equation}
    \gamma= \begin{cases}
        -f_{\tt acc}^l, & \mbox{if}~ -f_{\tt acc}^l \leq M \\
        M, &  -f_{\tt acc}^l > M ~\mbox{or}~ f(\bar{x}) - f(x^*) = 0,
    \end{cases}\label{eq:accuracymeasure}
    \end{equation}
or the analogous equation using $x_{\tt acc}^l$.  Note that we have multiplied $f_{\tt acc}^n$ by $-1$, so $\gamma$ can be interpreted as the number of new digits of accuracy obtained up to a maximal improvement value of $M$.

Similar measures can be used to quantify the amount of constraint violation for a test run.  Considering $\min \{ f(x) ~:~ g_i(x) \leq 0, i = 1, 2, ..., m\}$,
    $$\begin{array}{ll}
    \displaystyle\sum_{i=1}^m \max\{0, g_i(\bar{x})\} &\mbox{gives the sum of violated constraints},\\
    \displaystyle\sum_{i=1}^m (\max\{0, g_i(\bar{x})\})^2 &\mbox{gives the squared sum of violated constraints},\\
    \displaystyle\frac{1}{m} \sum_{i=1}^m \max\{0, g_i(\bar{x})\} &\mbox{gives the mean constraint violation, and}\\
    \displaystyle\prod_{i : g_{i}(\bar{x}) > 0} g_i(\bar{x}) &\mbox{amounts to the geometric mean of the violated constraints}.
    \end{array}$$

The selection of the appropriate strategy among the variety of approaches depends on the objectives of the experimental research, the problem structure, and the type of optimization algorithms used. The researcher should also carefully select the success criteria, e.g., how to fairly compare a solution that barely satisfies the constraints versus a solution that barely violates the constraints but returns a much lower objective function value.

\textbf{No known solution available:} In many situations, the test set used will not have known solutions to all problems.  This is particularly true if the test set includes instances of real-world applications.  To evaluate the quality of an algorithmic output in this situation, some new considerations are required \cite{114,115}.

First, it should be immediately obvious that, if no known solution is available, then fixed-target approaches cannot be applied.  Fixed-cost approaches are still applicable, but since $f(x^*)$ is not known, measuring the accuracy of the final algorithmic output's function value, $f(\bar{x})$, becomes difficult.  Measuring the accuracy of the final algorithmic output's point, $\bar{x}$, becomes essentially impossible.

To measure the quality of the final algorithmic output's function value $f(\bar{x})$, the simplest approach is to replace $f(x^*)$ with the best known value for the problem.  For any given test run, this guarantees that at least one algorithm gets the exact answer, so it is important to select a reasonable maximal improvement value.  Another approach is to estimate the optimal solution using statistical techniques. For example, in combinatorial optimization problems, some researchers \cite{c11, c12} use a sample of algorithmic outputs to predict the location of the solution. In \cite{c2}, such an approach is explained in an evaluation of non-deterministic algorithms.  Another strategy is to calculate a lower bound on the cost of an optimal solution, and to compare the algorithmic output cost with that lower bound. As an example, the total sum of the weight list in packing problems can be considered as a lower bound on the total number of bins used in a packing. Finally, one may abandon comparing the algorithmic output quality with the optimal solution, and only assess the quality of the algorithmic output with similar results published in the literature or other algorithms being tested.

\subsection{Parameter tuning and stopping conditions}

Additional parameters, such as stopping tolerances, population size, step sizes, or initial penalty parameters, are required for most optimization algorithms. 

Among such parameters, stopping conditions play a highly notable role, as different stopping conditions can drastically change the output of an algorithm \cite{SS00, SSL13, ZZ08}.  Moreover, if stopping tests are internalized within a method, it may not be possible to ensure all algorithms use the same stopping conditions \cite{SS00, SK06}.  However, if a fixed-cost or fixed-target approach (see Subsection \ref{ss:qualityofsolution}) is employed, then other stopping conditions can be turned off, thereby ensuring all algorithms use the same stopping conditions.   If it is not possible to ensure all algorithms use the same stopping conditions, then researchers should recognize this potential source of error when drawing conclusions from the results.

Other parameters, such as initial step length, can also have an impact on the performance of an optimization algorithm.  Such parameters often require tuning in order to obtain a better performance. If different choices of input parameters are allowed in an algorithm, researchers should mention the parameter settings used and how they were selected. Different strategies used for tuning parameters affect the benchmarking process. Choosing appropriate parameter settings for an optimization algorithm is usually based on experiments and statistical analysis.

The tuning strategy should be chosen in conjunction with a specific algorithm and in a replicable manner \cite{115}. Any improvements obtained from hand-tuning can of course be reported, but separately from more systematic comparative experiments. In some studies, algorithmic methods are presented to automate the tuning procedure of parameters \cite{AudetOrban2006, AudetDangOrgan2014, 145, 147, 121, 76}. The major disadvantage of these tuning methods is that they require a considerable computational investment because they usually try many possible settings to find an appropriate one. Nonetheless, in recent years some studies have specifically focused on automatic tuning of parameters in optimization solvers. Examples of these efforts include the machine learning based method proposed in \cite{145}, CPLEX automatic tuning tool \cite{c13}, use of derivative-free optimization \cite{AudetOrban2006}, ParamILS \cite{147}, and the procedure proposed in \cite{146} for mixed integer programming solvers. Similarly, some of the tuning techniques for non-deterministic methods include sequential parameter optimization (SPO) \cite{143, 113}, relevance and calibration approach \cite{144}, and F-Race \cite{c9}.

In view of the considerable research on automatic tuning of optimization solvers, a more accurate approach for benchmarking of optimization solvers requires a pre-processing step in which an automatic tuning method is employed to find the suitable configuration settings for all the solvers \cite{82}.  As this is not always practical, it is important to emphasize that tuning parameters can have a major impact on the performance of an algorithm, therefore it is not appropriate to tune the parameters of some methods while leaving other methods at their default settings.
%All algorithms should have parameters tuned in a similar manner and to equal degrees.

\section{Analyzing and reporting the results} \label{sec4}

Many studies use basic statistics (e.g., average solving time) to report the experimental results.  Basic statistics are a reasonable starting point, but provide little information about the overall performance of optimization methods.  Reporting methods can be loosely broken down into three categories: numerical tables, graphics, and performance ratio methods (e.g., performance and data profiles).

\subsection{Tables}

Numerical tables provide the most complete method of reporting benchmarking results, so for the sake of completeness, we recommend making full tables of results readily available.  However, such tables are cumbersome, so are often better included in an appendix or in additional online material linked to an article.

As full tables of results can easily overwhelm a reader, researchers have developed various techniques that provide easy-to-understand and compact methods for reporting the experimental results.  Summary tables can be employed as a first step \cite{SK06}.  For example, in \cite{5} optimization methods were rated by the percentage of problems for which a method's time is termed \textit{competitive} or \textit{very competitive}.  The solving time of an algorithm was called competitive if $T_s\le 2 T_{min}$ in which $T_s$ is the solving time obtained by an algorithm on a particular problem and $T_{min}$ is the minimum solving time obtained among all the algorithms on that specific problem.  Similarly, if $T_s\le \frac{4}{3} T_{min}$, then they call that method very competitive.  Tables such as these provide good talking points for discussing benchmarking data, but fail to give a complete picture of the results.  One critic for this particular approach is it does not explore how much the table would change if, for example, the cut-off for very competitive was changed from $\frac{4}{3} T_{min}$ to $\frac{5}{4} T_{min}$.

Many other forms of summary tables are present throughout optimization benchmarking literature, however all suffer from the same fundamental problem -- to be readable a summary table must distill the results down to a highly condensed format, thereby eliminating much of the benchmarking information.

\subsection{Graphics}

Well-conceived graphics can provide more information than some other data presentations. Simple graphical methods, such as histograms, box-plots, and trajectory plots, provide a next step in analysis, while more complete methods include performance profiles, data profiles, and accuracy profiles. Depending on the objectives of an experimental research, one or more of these techniques might be useful to report the results. In \cite{132, 133}, different types of plots are introduced, which are useful for data representation in general.

\begin{figure}[ht]
\captionof{figure}{A sample trajectory plot.}
\centering
\fbox{\includegraphics[width=13cm]{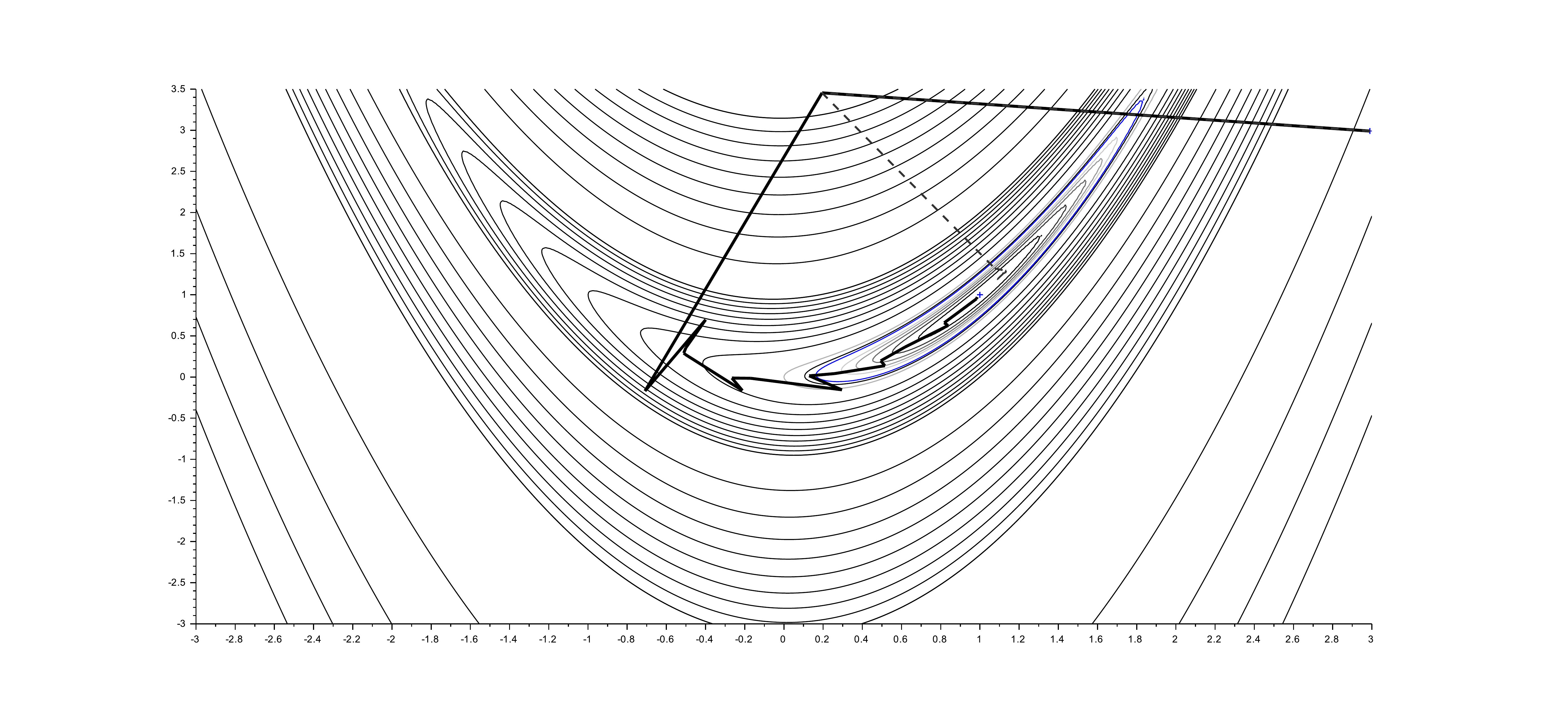}}
\label{Fig00}
\end{figure}

A more specialized plot for optimization algorithms is the trajectory plot \cite{1, 22, 108, 105, 106, c8, SS00}.  In a trajectory plot, the performance of an optimization algorithm on a given test problem is visualized by plotting a path that connects the points generated by each iteration of the algorithm.  An example appears in Figure \ref{Fig00}, where the trajectories of two algorithms attempting to minimize the Rosenbrock function are plotted.  Both algorithms begin at the point $(3, 3)$, and the first iteration moves both algorithms to the point $(0.2, 3.5)$.  Algorithm 1 (represented by the solid line) proceeds to $(0.7, -0.2)$ and continues in a zig-zag path towards the minimizer.  Algorithm 2 (represented by the dashed line) proceeds to $(1.1, 1.3)$ and then follows a fairly straight path towards the minimizer, albeit with very small step sizes.  While trajectory plots are useful to build a better understand of how each algorithm behaves, they are not particularly good for benchmarking as they can only present the results for one test problem at a time. They are also limited to plots of functions of 2 or 3 variables, or to plotting projections onto subspaces for more than 3 variables.

\begin {figure}[H]
\captionof{figure}{A sample convergence plot.}
\centering
\fbox{\includegraphics[scale=0.5]{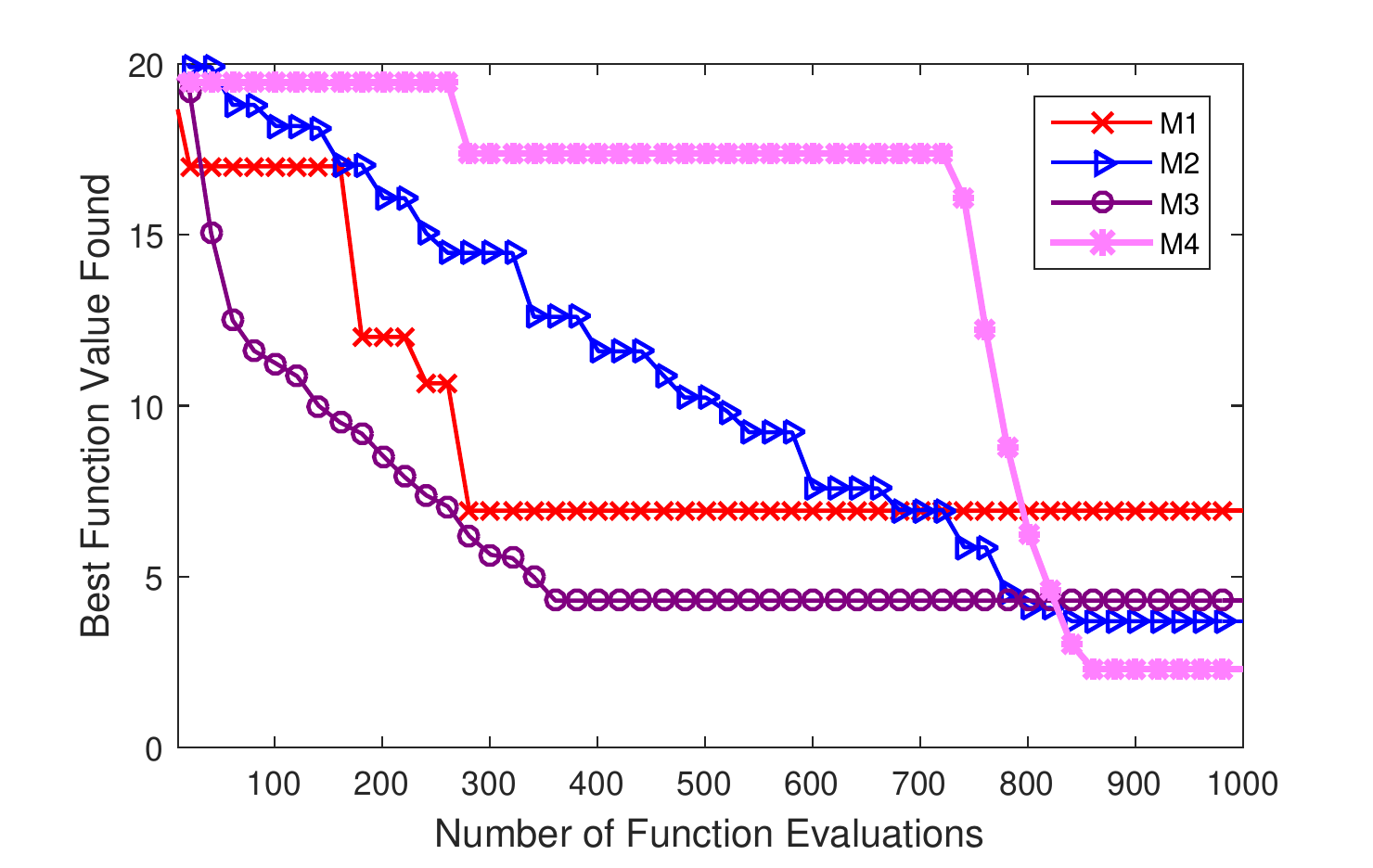}}
\label{Fig01}
\end {figure}

Another specialized plot for optimization benchmarking is the {\em convergence plot}.  In a convergence plot the performance of different optimization methods is visualized by plotting the best function value found against some measure of fundamental evaluation (Section \ref{ss:efficiency}).  An example convergence plot is given in Figure \ref{Fig01}.

In Figure \ref{Fig01} the results of four optimization methods are plotted for a given test problem.  In this example, method M1 starts well, but stalls after about 300 function evaluations, while method M2 shows steady decrease for about 800 function evaluations before stalling.  Method M3 initially decreases the fastest, but stalls after about 350 function evaluations. Finally, method M4 starts very slowly, but ultimately finds the lowest value.  Like trajectory plots, convergence plots are useful for discussing some specific behavior of the algorithm, but are poor for benchmarking as they can only be used to analyze one test problem at a time.

While trajectory and convergence plots are useful to visualize a method on one problem, their main drawback is that they represent the results for a single problem per plot.  So if the test set contains a large number of problems then it will be difficult to evaluate the overall performance of these methods. Other types of plots can be found in the literature, but generally have the same limitations as trajectory and convergence plots \cite{17} \cite{22}.

For many optimization algorithms, researchers are interested in how the problem scales with the size of the input (e.g., dimension of the problem). For such research it can be valuable to produce a {\em runtime plot}. Runtime plots visualize the data by plotting the time to solve across a series of problem instances with different sizes.   Runtime plots can suffer from a similar issues to trajectory and convergence plots, namely, they represent the results for a single series of problem instances.  However, this problem can be somewhat mitigated by aggregating data from a collection of problems to create an ``average runtime'' plot.

\subsection{Performance profiles}

According to \cite{SKM16}, the idea of creating graphical comparison of optimization methods dates back to at least 1978 in the paper by Grishagin \cite{G78}.\footnote{We thank ``Mathematics Referee \#1'' for pointing out that reference.}  In 2000, Strongin and Sergeyev presented the idea of {\em operational characteristics} for an algorithm: a graphical method to visualize the probability that an algorithm solves a problem within a set time-frame \cite{SS00}.  However, it was not until the 2002 paper by Dolan and Mor\'e, \cite{2}, that the idea of graphically presenting benchmarking results became mainstream.  Dolan and Mor\'e (apparently unaware of works of Grishagin or Strongin and Sergeyev) denoted their proposed graphs {\em performance profiles}.

Performance profiles provide interesting information such as efficiency, robustness, and probability of success in a graphically compact form \cite{2}.  Their use has grown rapidly in optimization benchmarking, and should certainly be considered for any benchmarking optimization research.  

Let $\mathcal{P}$ be a set of problems, $\mathcal{S}$ a set of optimization solvers, and $\mathcal{T}$ a convergence test. Assume proper data has been collected. Performance profiles are now defined in terms of a performance measure $t_{p,s} > 0$, obtained for each pair of $(p,s)\in P\times S$.  This measure can be the computational time, the number of function evaluations, etc. A larger value of $t_{p,s}$ shows worse performance. For each problem $p$ and solver $s$, the performance ratio is defined as

\begin{equation}\label{05}
  r_{p,s}=
  \begin{cases}
    \displaystyle \frac{t_{p,s}}{\min\{t_{p,s}:s \in \mathcal{S}\}} &\text{if  convergence test passed,} \\
    \displaystyle \infty &\text{if convergence test failed.}
  \end{cases}
\end{equation}
for a specific problem $p$ and test $s$ (the best solver has $r_{p,s} = 1$). The \textit{performance profile} of a solver $s$ is defined as follows

\begin{equation} \label{06}
\rho_s({\tau})=\frac{1}{\lvert \mathcal{P} \rvert} \text {size} \{\mathit{p}  \in \mathcal{P}: r_{p,s} \le \tau \},
\end{equation}
where $\lvert \mathcal{P} \rvert$ represents the cardinality of the test set $\mathcal{P}$. Then, $\rho_s({\tau})$ is the portion of the time that the performance ratio $r_{p,s}$ for solver $s \in \mathcal{S}$ is within a factor $\tau \in \mathbb{R}$ of the best possible performance ratio.

Note that $\rho_s({1})$ represents the percentage of problems for which solver $s \in \mathcal{S}$ has the best performance among all the other solvers. And for $\tau$ sufficiently large, $\rho_s({\tau})$ is the percentage of the test set that can be solved by $s \in \mathcal{S}$. Solvers with consistently high values for $\rho_s({\tau})$ are of interest.

\begin {figure}[H]
\caption{An example performance profile using logarithmic values of $\tau$.}
\centering
\fbox{\includegraphics[scale=0.5]{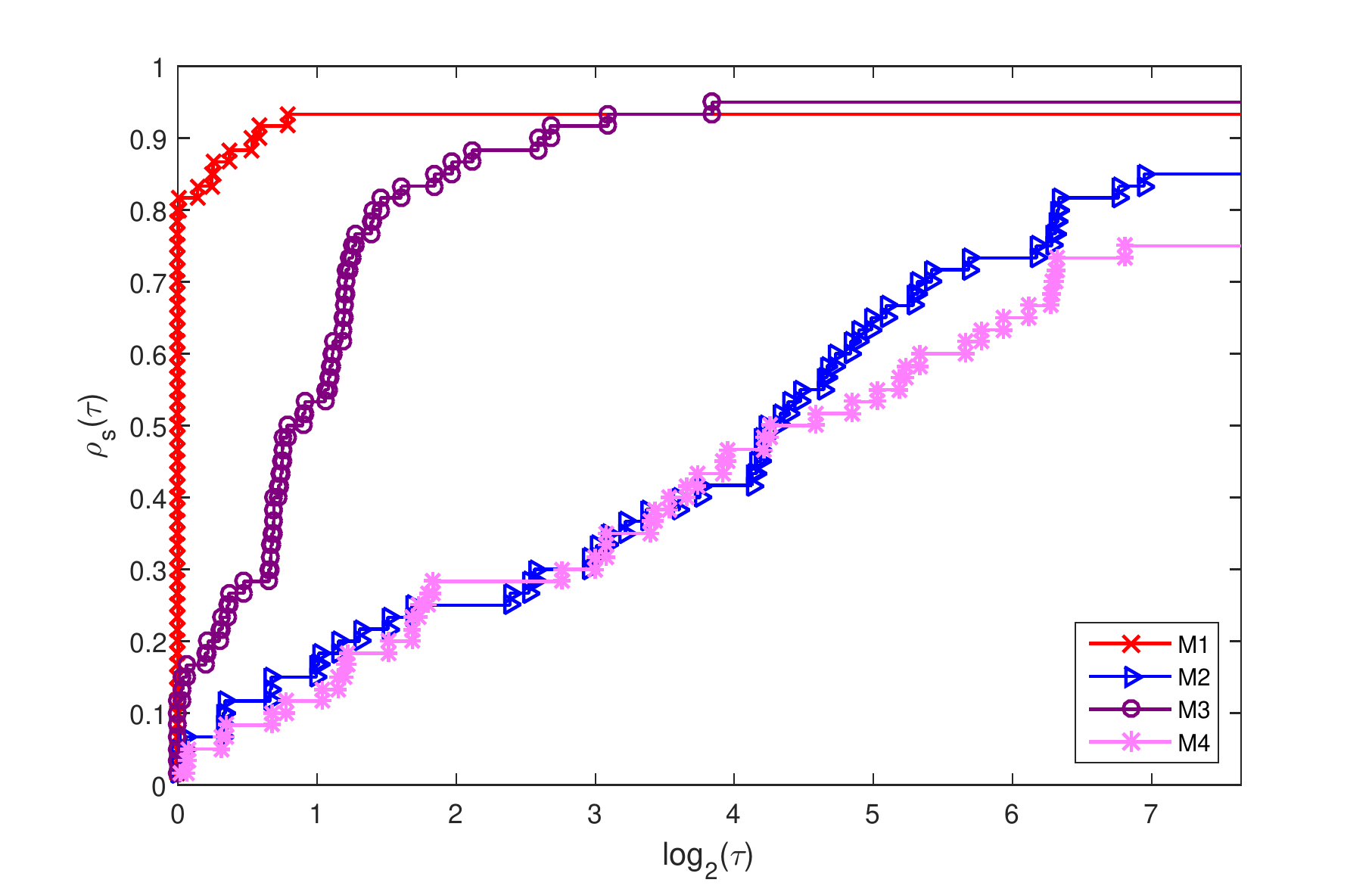}}
\label{Fig02}
\end {figure}

\begin {figure}[H]
\caption{The performance profile from Figure \ref{Fig02} using non-logarithmic values of $\tau$.}
\centering
\fbox{\includegraphics[scale=0.5]{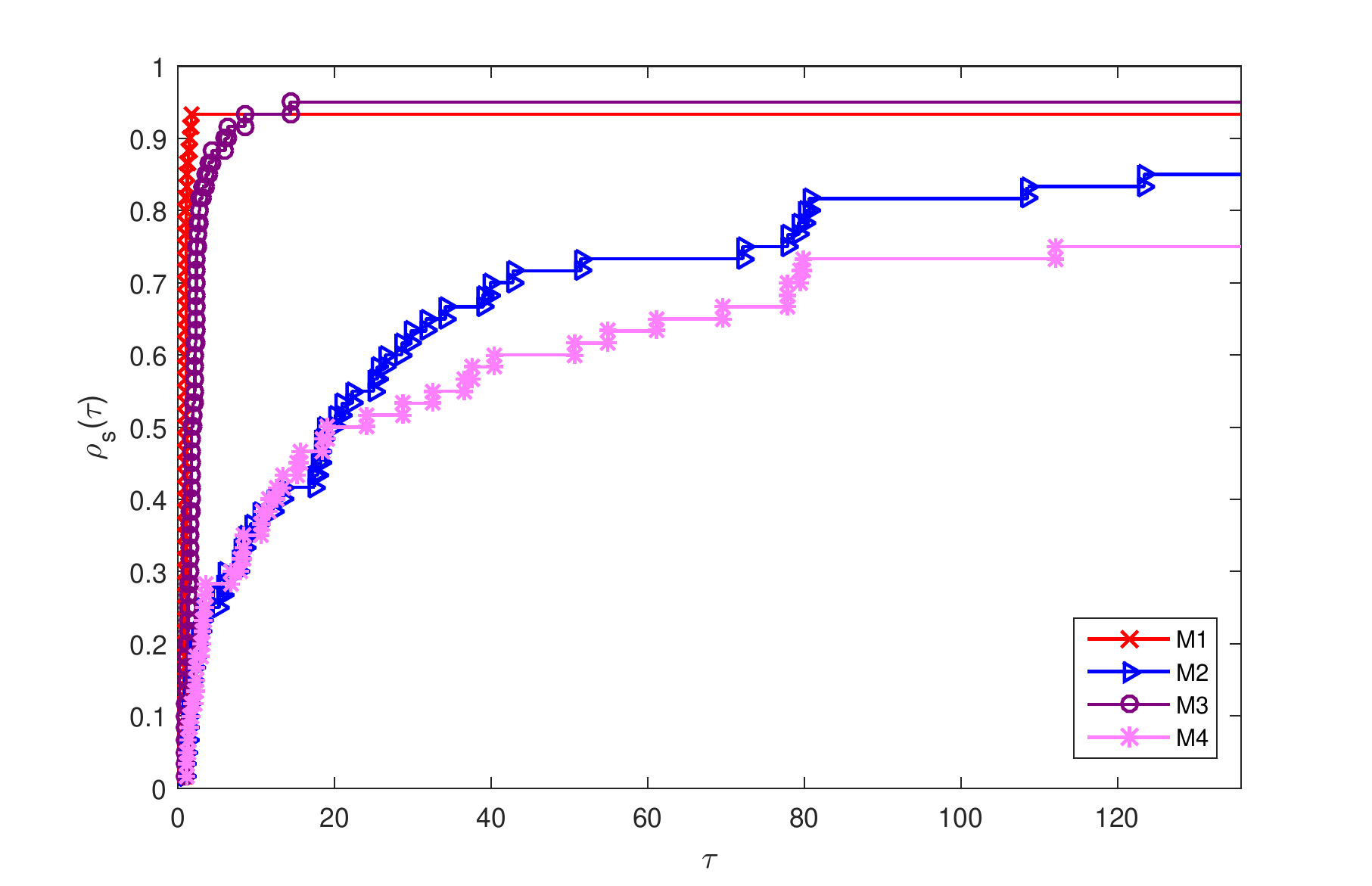}}
\label{Fig02-2}
\end {figure}

Figure \ref{Fig02} shows a sample performance profile plot (created using data from \cite{c7}) for logarithmic values of $\tau$. The logarithmic values are employed to deal with smaller values for $\tau$. This will result in a more accurate plot which shows the long term behavior of the methods. To demonstrate the difference, Figure \ref{Fig02-2} shows the same performance profile using non-logarithmic values of $\tau$.  Depending on the data collected, logarithmic or non-logarithmic values of $\tau$ may be more appropriate.  Researchers should create both profiles, but it may be only necessary to provide one in the final manuscript.

The performance profiles in Figure \ref{Fig02} compares four different optimization methods on a test set of 60 problems. The method M1 has the best performance (in terms of CPU time) for almost 93\% of the problems; meaning that M1 is able to solve 93\% of the problems as fast or faster than the other two approaches. M3 solves roughly 11\% of the problems as fast or faster than the other approaches. On the other hand, given enough time M1 solves about 92\% of all problems, while M3 solves about 94\% of all problems.  The graphs of M1 and M3 cross at about $\log_2(\tau) \approx 3$ (i.e., $\tau \approx 8$), the two methods solve the same number of problems if time to solve is relaxed to be within a factor of $8$. 
%After that point, M2 and M4 have a similar performance up to $\log_2(\tau)=4$. After that point, M2 has a better performance than M4 and ends up solving 85\% of the problems which is around 10\% more than M4.

Since performance profiles compare different methods versus the best method, the interpretation of the results should be limited to comparison to the best method and no interpretation should be made between, e.g., the second best and third best method since a switching phenomenon may occur~\cite{GOULD-16}\footnote{We thank ``Engineer Referee \#3'' for pointing out that reference.}. To compare the second and third best method, a new performance profile should be drawn without the first method, see the explicit examples provided in~\cite{GOULD-16}.

Performance profile plots can be customized by substituting the standard performance measure \textit{time}. For example, in \cite{4, 134, 136}, the objective function value is used as the performance measure to compare the profiles. In particular, $t_{p,s}$ is replaced with 
\begin{equation}
m_{p,s} = \frac{\hat{f}_{p,s}(\text{after $k$ function evaluations}) - f^*}{(f_w-f^*)},
\end{equation}
for problem $p$ and solver $s$, where $f_w$ is the largest (worst) function value obtained among all the algorithms, and $\hat{f}_{p,s}$ is the estimated function value after $k$ function evaluations.  In another example, \cite{SK15} creates a performance measure based on proximity to optimal points.

The primary advantage of performance profiles is that they implicitly include both speed and success rate in the analysis.  The value of $\rho_s({\alpha})$ gives a sense of how promising the algorithmic outputs are relative to the best solution found by all the optimization algorithms that are compared together.  

One criticism of performance profiles is that the researcher must select a definition for the convergence test passing and failing.  Changing this definition can substantially change the performance profile \cite{HareKochLucet2011}.  Also note that if a fixed-cost approach is used to performing the benchmarking experiments, then performance profiles become inappropriate, as all algorithms will use the same ``time''. Another criticism is that the profile is only showing performance with respect to the best method and does not allow one to compare other methods with each other due to the appearance of a switching phenomenon~\cite{GOULD-16}. 

Nonetheless, performance profiles have become a gold-standard in modern optimization benchmarking, and should be included in optimization benchmarking analysis whenever possible with an appropriate interpretation.

\subsection{Accuracy profiles}

\begin {figure}[H]
\caption{An example accuracy profile.}
\centering
\fbox{\includegraphics[scale=0.5]{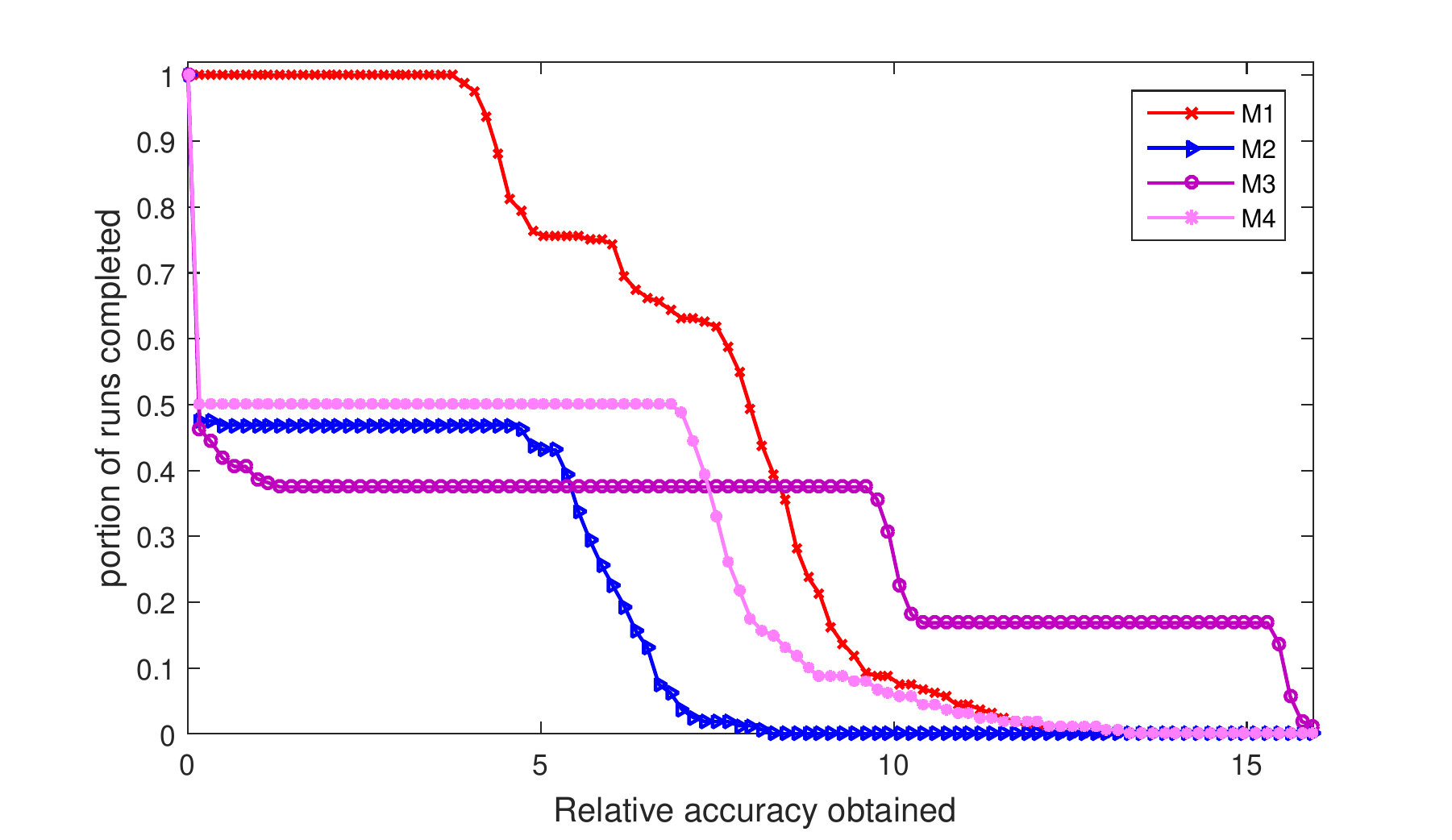}}
\label{Fig05}
\end {figure}

Similar to performance profiles, \textit{accuracy profiles} provide a visualization of an entire optimization benchmarking test set.  However, accuracy profiles are designed for fixed-cost data sets.  They begin by defining, for each problem $p \in \mathcal{P}$ and solver $s \in \mathcal{S}$, an accuracy measure (similar to equation \eqref{eq:accuracymeasure}):
   $$\gamma_{p, s} = \begin{cases}
        -f_{\tt acc}^{p,s}, & \mbox{if}~ -f_{\tt acc}^{p,s} \leq M \\
        M, &  -f_{\tt acc}^{p,s} > M ~\mbox{or}~ f_{\tt acc}^{p,s} ~\mbox{is undefined},
    \end{cases}$$
where $f_{\tt acc}^{p,s} = \log_{10}(f(\bar{x}_{p,s}) - f(x^*_p)) - \log_{10}(f(x^0_p) - f(x^*_p))$, $\bar{x}_{p,s}$ is the candidate solution point obtained by solver $s$ on problem $p$, $x^*_p$ is the optimal point for problem $p$, and $x^0_p$ is the initial point for problem $p$.  The performance of the solver $s \in \mathcal{S}$ on the test set $\mathcal{P}$ is measured using the following function
    $$
    R_s (\tau)= \frac{1}{|\mathcal{P}|} \text {size} \{\gamma_{p,s} | \gamma_{p,s} \ge \tau, p \in \mathcal{P} \}.
    $$
The accuracy profile $R_{s} (\tau)$ shows the proportion of problems such that the solver $s \in \mathcal{S}$ is able to obtain a solution within an accuracy of $\tau$ of the best solution. An example accuracy profiles (using data from \cite{Hare-Sagastizabal-2006}) appears in Figure \ref{Fig05}.

In Figure \ref{Fig05}, we see 4 methods (M1, M2, M3, and M4) plotted against each other in an accuracy profile format.  Examining the profile, notice that method M1 achieves 5 digits of accuracy on almost all test-problems, and 6 digits of accuracy on about 75\% of test problems.  All other method achieve this level of accuracy on 50\% or less of test problems.  Thus, if 5 or 6 digits is the desired level of accuracy, then M1 is a clear winner.  However, if the particular application requires much higher accuracy, then M3 becomes a contender.  Indeed, only M3 was able to achieve 12 digits of accuracy on any reasonable portion of the test problems.  (In this particular test, accuracy was capped at 16 digits, but no method managed to achieve this on a significant portion of the test problems.)

Accuracy profiles do not provide as much information as performance profiles, but are suitable when fixed-cost data sets are collected.  This is appropriate in cases where the cost of obtaining the exact solution exceeds the budget, so the optimization target is to find as good a solution as possible within a limited time.

\subsection{Data profiles}

%In another approach to demonstrate the performance of optimization algorithms in a compact and easy-to-understand form, 
Mor\'{e} and Wild \cite{3} proposed data profiles as an adjustment to performance profiles for derivative-free optimization algorithms. Data profiles try to answer the question: what percentage of problems (for a given tolerance $\tau$) can be solved within the budget of $k$ function evaluations? They assume the required number of function evaluations to satisfy the convergence test is likely to grow as the number of variables increases. The data profile of an optimization algorithm $s$ is defined using~\cite{3}
\begin{equation} \label{07}
d_s({k})=\frac{1}{\lvert \mathcal{P} \rvert} \text {size} \left\{\mathit{p}  \in \mathcal{P}: \frac{t_{p,s}}{n_p+1} \le k \right\},
\end{equation}
in which $t_{p,s}$ shows the number of function evaluations required to satisfy the convergence test, $n_p$ is the number of variables in the problem $p \in \mathcal{P}$, and $d_s(k)$ is the percentage of problems that can be solved with $k(n_p + 1)$ function evaluations.  The value $k(n_p + 1)$ is used since $n_p + 1$ is to the number of function evaluations required to compute a ``simplex gradient'' (a one-sided finite-difference estimate of the gradient).  

It is worth noting that data profiles could easily be defined replacing $\frac{t_{p,s}}{n_p+1}$ by any other measure of fundamental evaluations used.  Moreover, if $\frac{t_{p,s}}{n_p+1}$ is replaced by iterations, then data profiles become a slight variation of {\em operational characteristics} defined in \cite{SS00}.

\begin {figure}[H]
\caption{An example data profile.}
\centering
\fbox{\includegraphics[scale=0.5]{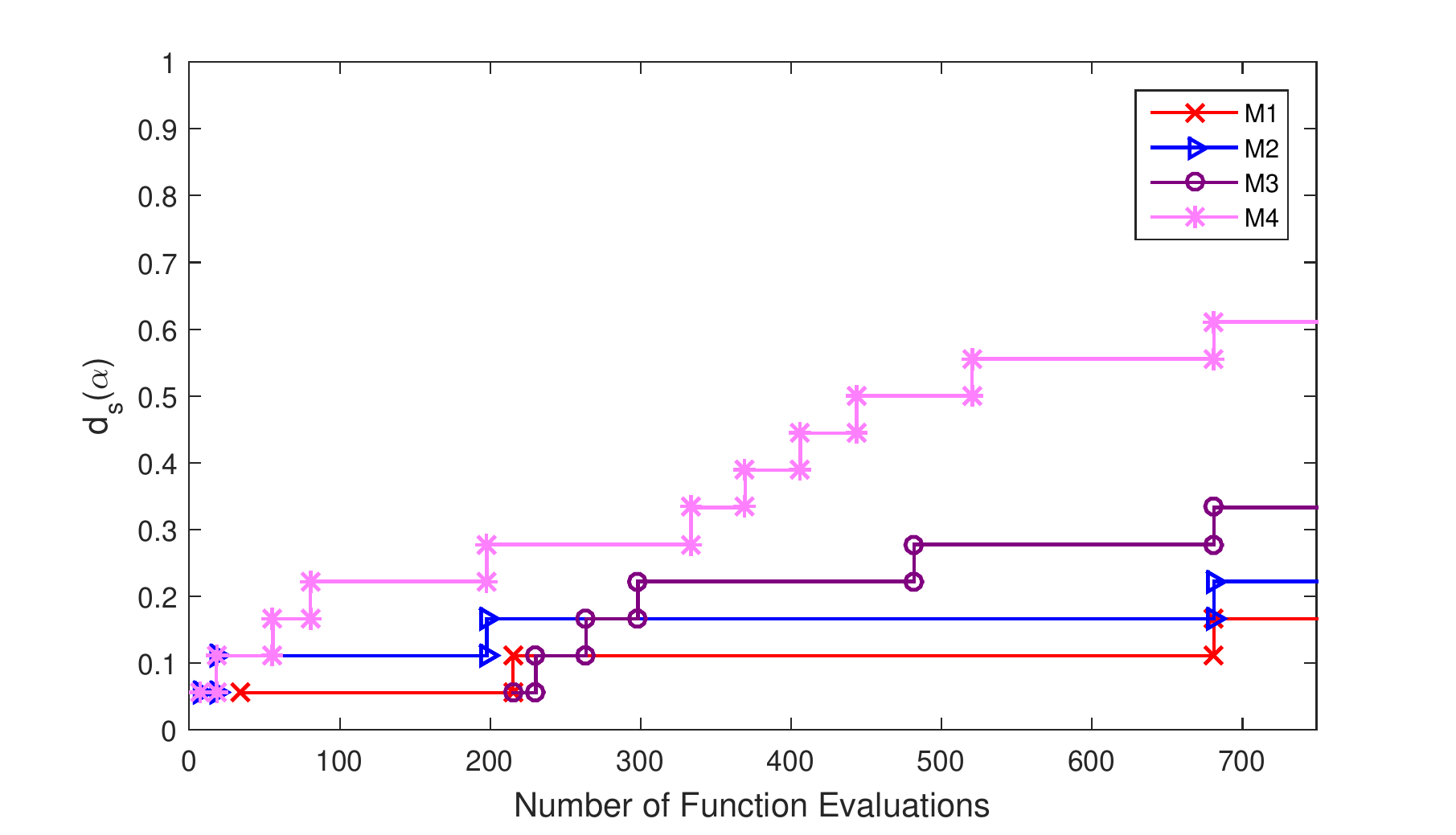}}
\label{Fig03}
\end {figure}

Figure \ref{Fig03} shows a typical data profile. Suppose the user has a budget limit of 100 simplex gradients, according to Figure~\ref{Fig03}, with this budget the method M4 has the best performance by solving roughly $22\%$ of the problems; while M3 has the worst performance among all the solvers since with this budget it does not solve any problem.

Like performance profiles, data profiles are cumulative distribution functions, and thus, monotone increasing step functions with a range in $[0, 1]$.  Data profiles do not provide the number of function evaluations required to solve a specific problem, but instead provide a visualization of the aggregate data.  Also note that the data profile for a given solver $s \in S$ is independent of other solvers; this is not the case for performance profiles.

Although the data profiles are useful for benchmarking, they have not received the same extensive attention as the performance profiles.  This is partly because they are newer, but perhaps also because they are primarily used with derivative-free optimization algorithms.  However, data profiles could be easily adjusted to a broader class of algorithms by replacing $t_{p,s}$ with any measure of time, and $n_p + 1$ by any dimensional normalization factor.  For example, for a sub-gradient based method, $d_s(\alpha)$ could be redefined as
    $$d_s({\alpha})=\frac{1}{\lvert \mathcal{P} \rvert} \text {size} \{\mathit{p}  \in \mathcal{P}: g_{p,s} \le \alpha \},$$
where $g_{p,s}$ is the number of sub-gradient evaluations.  This might make them an appropriate tool for benchmarking bundle methods \cite[\S XIV-XV]{HUL93}.

Table \ref{Table4-Report} summarizes the reporting methods discussed in this section.

\begin {sidewaystable}[H]
\captionof{table}{Reporting methods summarization.\label{Table4-Report}}
\begin{center}
\begin{tabular*}{\textwidth}{lllll}
  \toprule
	\head{Reporting method} & \head{Evaluates} & \head{Advantage} & \head{Drawback} & \head{Recommendation} \\
  \midrule
      \begin{tabular}{m{3.0cm}}Full Data Tables\-\end{tabular}
    & \begin{tabular}{m{3.5cm}}---\end{tabular}
    & \begin{tabular}{m{3.8cm}}Comprehensive \end{tabular}
    & \begin{tabular}{m{3.5cm}}Overwhelming \end{tabular}
    & \begin{tabular}{m{4.0cm}}Provide in appendix or online data set.\end{tabular} \\
  \midrule
      \begin{tabular}{m{3.0cm}}Summary Tables \\\\ Simple Graphs \end{tabular}
    & \begin{tabular}{m{3.5cm}}Varies\end{tabular}
    & \begin{tabular}{m{3.8cm}}Brief\end{tabular}
    & \begin{tabular}{m{3.5cm}}Incomplete\end{tabular}
    & \begin{tabular}{m{4.0cm}}Provide as talking point, \\ but include other forms of analysis.\end{tabular} \\
  \midrule
      \begin{tabular}{m{3.0cm}} Trajectory Plots \\\\ Convergence Plots\end{tabular}
    & \begin{tabular}{m{3.5cm}} Speed and Accuracy \\\\ Efficiency \end{tabular}
    & \begin{tabular}{m{3.8cm}} Clear \\\\ Precise \end{tabular}
    & \begin{tabular}{m{3.5cm}} Examines one problem at a time.\end{tabular}
    & \begin{tabular}{m{4.0cm}} Good for case-studies, but should include other forms of analysis for larger data sets\end{tabular} \\
  \midrule
      \begin{tabular}{m{3.5cm}} Performance Profiles\- \\\\ Accuracy Profiles\- \\\\ Data Profiles\end{tabular}
    & \begin{tabular}{m{3.5cm}} Speed and Robustness \\\\Accuracy\\\\Speed and Robustness\end{tabular}
    & \begin{tabular}{m{3.8cm}} Strong graphical representation that incorporates the entire dataset.\end{tabular}
    & \begin{tabular}{m{3.5cm}} Cannot be used for fixed-cost data sets\end{tabular}
    & \begin{tabular}{m{4.0cm}} Include at least one of these three profiles whenever possible\end{tabular}  \\
  \bottomrule
\end{tabular*}
\end{center}
\end{sidewaystable}

\section{Automated benchmarking}\label{sec5}

As we have seen, the benchmarking process of optimization algorithms is a complicated task that requires a lot of effort from data preparation and transformation to analysis and visualization of benchmarking data. Accordingly, some researchers have begun the development of software tools to facilitate and automate developing test sets, solving the problems using a variety of optimization algorithms, and carrying out performance analysis and visualization of benchmarking data.

The PAVER server \cite{33,45} is an online server that provides some tools for automated performance analysis, visualization, and processing of benchmarking data. An optimization engine, either a modeling environment, such as AMPL \cite{127} or GAMS \cite{126}, or a stand-alone solver, is required to obtain solution information such as objective function value, resource time, number of iterations, and the solver status. Then, the benchmark data obtained by running several solvers over a set of problems can be automatically analyzed via online submission to the PAVER server. PAVER returns a performance analysis report through e-mail in HTML format. The tools available in PAVER allow either direct comparisons between two solvers or comparisons of more than two solvers simultaneously in terms of efficiency, robustness, algorithmic output quality, or performance profiles.

The High-performance Algorithm Laboratory~\cite{44} (HAL) is a computational environment designed to facilitate empirical analysis and design of algorithms. It supports conducting large computational experiments and uses a database to handle data related to algorithms, test sets, and experimental results. It also supports distributed computation on a cluster of computers. Its major advantage over other tools is its aim to develop a general purpose tool that can handle different categories of problems, although the initial deployment of problems and algorithms is tricky.

The Optimization Test Environment \cite{74} is another tool that can be used for benchmarking the performance of different optimization algorithms. It provides some facilities to organize and solve large test sets, extract a specific subset of test sets using predefined measures, and perform statistical analysis on the benchmarking data. The results obtained by each optimization algorithm is verified in terms of feasibility and correctness. A variety of information is reported such as the number of global numerical solutions found (i.e., the best solution found among all optimization algorithms), number of local solutions found, number of wrong claims, etc. For problem representation, it uses Directed Acyclic Graphs (DAGs) from the Coconut Environment \cite{128}. This user-friendly environment analyzes results and automatically summarizes them before reporting them in an easy-to-use format such as \LaTeX, JPEG, and PDF.

Other software tools for automating benchmarking process include EDACC \cite{140}, LIBOPT \cite{46}, CUTEr \cite{87} and a testing environment reported in \cite{5}.

Using automated performance analysis tools has the potential to facilitate the benchmarking process.  Moreover, the automation of the process may reduce the risk of biased comparison, by taking some of the comparison decisions away from the algorithm designer.  However, automated benchmarking tools are not yet accepted by the research community due to their shortcomings. The major drawback of these tools is that the flexibility of a researcher to design experiments based on their research objectives is restricted to the tools' limitations and the way they view the benchmarking process. Moreover, so far all of these tools operate in expert mode, meaning that the usability aspect needs to be improved in terms of application and design of experiments. In most cases preparation of an experiment beyond the scope of default facilities of the benchmarking tools is nontrivial and involves some customization, e.g., scripting. Further research in this direction will create valuable tools for the optimization community, but the current status is not ready for wide-spread use.

\section{Conclusion} \label{sec6}

This article reviews the issue of benchmarking optimization algorithms. For the sake of having a careful, less-biased, explicitly-stated, and comprehensive evaluation of the optimization algorithms an {\em a priori} benchmarking design is required. To develop an appropriate experimental design, the first thing to consider is to clarify the questions that are to be answered by the experiment. This includes selecting a suitable test set and suitable performance measures based on the objectives of the research. The data must be analyzed and processed in a transparent, fair, and complete manner.  Within this paper we discuss each of these topics, and present a review of the state-of-the-art for each of these steps.  We include several tables and figures that summarize the process, and provide key advice designed to lead to a fair benchmarking process.

A final important point must be raised in regards to optimization benchmarking:
    \begin{quote}{\bf as in all scientific research, benchmarking optimization algorithms should be reported in a manner that allows for reproducibility of the experiments.}\end{quote}
When reporting results, make sure to describe algorithms, parameters, test problems, the computational environment, and the statistical techniques employed with an acceptable level of details.  It should be clarified that it is usually difficult to provide enough information in a published paper to enable the reader to rerun the stated experiments and replicate completely the reported results. Moreover, the pace of computational development is so high that it is virtually impossible to entirely reproduce a computational experiment, due to development and modifications in operating systems, computer architecture, programming languages, etc.  However, the minimum standard for replication of the experiments is that at least the authors themselves should be able to replicate the experiments \cite{38}.  Therefore, it is important that the researcher  keep all the programs and data necessary to redo all the computations and recreate all graphs.  Such programs should be made available whenever possible.

\subsection{Some final insights and remarks from the referees}

This paper provides a high-level perspective on benchmarking of optimization algorithms.  While it does not aim to be all encompassing, it hopefully provides a baseline for best practices when benchmarking optimization algorithms.  Many nuances exist when dealing with benchmarking specific genres of optimization algorithms. We end with some final discussion on some of these nuanced areas.  Many of these final remarks were provided through the insights of $5$ excellent referees.

The state-of-the-art in optimization benchmarking currently has (at least) two major voids that require further research: how to properly benchmark optimization algorithms that make use of parallel processing, and how to properly benchmark multi-objective optimization algorithms.

Evaluating the performance of parallel optimization algorithms is different from traditional optimization methods in various aspects: performance measures such as time, the appropriate test sets, the new measures of merit involved in parallel processing such as the concept of speedup, efficiency, etc. All of these concerns together with the fast pace of technological advances in parallel computing motivate a research on the benchmarking of parallel optimization algorithms. A good start in this regard is the research paper by Barr and Hickman \cite{119}.

Benchmarking multi-objective optimization algorithms is similarly in its infancy. Appropriate test sets and performance measures have yet to surface.  Multi-objective optimization is a rapidly advancing field, and research into proper benchmarking in this discipline would be highly valuable.

A benchmarking challenge that we have not addressed is how to compare optimization algorithms that are different in nature\footnote{We thank ``Mathematics Referee \#1'' for pointing out this challenge.}.  For example, consider the comparison of a deterministic and a non-deterministic method \cite{GK17,KM17}.  If the multiple repeats of the non-deterministic method are considered, is it fair to compare the average quality to the single run of the deterministic method.  Some ideas on this, including a proposed method for comparing deterministic and a non-deterministic methods, can be found in \cite{SKM16}.

Another benchmarking challenge that has not been fully address is how to compare algorithms that approach the same problem from fundamentally different view points\footnote{We thank ``Engineering Referee \#3'' and ``Mathematics Referee \#2'' for pointing out this challenge.}.  For example, when working with constrained optimization problems, some researchers have explored {\em infeasible point methods} while others have focused on {\em interior point methods}.  Infeasible point methods typically take a two phase approach, where one phase aims for decrease in function value and the second phase aims to improve feasibility.  Interior point methods assume a strictly feasible starting point and use some form of penalty function to maintain feasibility of all trial points.  Comparing these two styles of algorithms is very challenging, and possibly meaningless, as one assumes an infeasible starting point and the other assumes a feasible starting point. Other algorithms adopt a hybrid approach by approximating the feasible set with some tolerance~\cite{REGIS-17}; in that case, the tolerance parameter could greatly influence the result of the comparison.

A source of debate in benchmarking global optimization algorithms is how to deal with {\em rescaling of the domain}\footnote{We thank ``Mathematics Referee \#2'' for pointing out this challenge.}.  Many global optimization algorithms are designed with the baseline assumption that the optimization problem's constrained region is the unit hypercube $[0,1]^n$.  Of course, in practical applications this is not always true.  Some algorithms deal with this at the solver level, using the constraint set's diameter to select parameters like initial step lengths; while other algorithms deal with this at the problem level, assuming that the end-user will rescale the constraint set to be the unit hypercube (which is not always easy to do).  Comparisons of algorithms that place fundamentally different assumptions on the problem structure may impact the selection of an appropriate test set and may limit the conclusions one can draw from the numerical results.

Another potential limitations on what conclusions can be drawn from a numerical study is the sensitivity analysis of the parameters\footnote{We thank ``Mathematics Referee \#1'' for pointing out this challenge.}. A robust study should investigate a range of parameters and report on their impact on the validity of the conclusions. We leave the complexity of how best to report such information to future research.

\bibliography{Benchmarking}
\bibliographystyle{alpha}
\end {document}